\documentclass{amsart}

\newtheorem{thm}{Theorem}

\newtheorem{prop}[thm]{Proposition}

\newcommand{\ABS}[1]{{\left| #1 \right|}} % |1|
\newcommand{\PAR}[1]{{\left(#1\right)}} % (1)
\newcommand{\SBRA}[1]{{\left[#1\right]}} % [1]
\newcommand{\BRA}[1]{{\left\{#1\right\}}} % {1}
\newcommand{\bs}{\ensuremath{\backslash}} % \
\newcommand{\defeq}{:=}
\newcommand{\p}[4]{{#3}\!\left#1{#4}\right#2} 

\newcommand{\moyf}[1]{\mathbf{E}_{#1}}

\newcommand{\moy}[2]{\p(){\moyf{#1}}{#2}}
\newcommand{\entf}[1]{{\mathbf{Ent}}_{#1}}
\newcommand{\ent}[2]{\p(){\entf{#1}}{#2}}
\newcommand{\eentf}{\mathbf{N}}
\newcommand{\eent}[1]{\p(){\eentf}{#1}}
\newcommand{\sentf}{\mathbf{H}}
\newcommand{\sent}[1]{\p(){\sentf}{#1}}
\newcommand{\fishf}{\mathrm{J}}
\newcommand{\fish}[1]{\p(){\fishf}{#1}}
\newcommand{\fishmf}{{\mathbf{J}}}
\newcommand{\fishm}[1]{\p(){\fishmf}{#1}}
\newcommand{\Kf}{\mathbf{K}}
\newcommand{\K}[1]{\p(){\Kf}{#1}}
\newcommand{\ds}[1]{\ensuremath{\displaystyle{#1}}}
\newcommand{\isopf}{\mathbf{I}}
\newcommand{\isop}[1]{\p(){\isopf}{#1}}
\newcommand{\g}{{\gamma}}
\newcommand{\la}{{\lambda}}
\newcommand{\gn}{\g_n}
\newcommand{\lan}{\la_n}

\newcommand{\TR}{\mathrm{Tr\,}}

\newcommand{\LA}{\mathrm{\Delta}}
\newcommand{\GR}{\mathrm{\nabla}}
\newcommand{\SGf}[1]{\mathbf{P}_{#1}}
\newcommand{\SG}[2]{\p(){\SGf{\!#1}}{#2}}

\newcommand{\dR}{\ensuremath{\mathbf{R}}}
\begin{document}

\title{Gaussian maximum of entropy and \\ reversed log-Sobolev inequality}
\author{Djalil Chafa\"\i}
\email{chafai@cict.fr}
\address{Laboratoire de Statistique et Probabilit\'es,%
  UMR CNRS 55830 \\ Universit\'e Paul Sabatier, F-31062 CEDEX, Toulouse, France.}

\begin{abstract}
  The aim of this note is to connect a reversed form of the Gross logarithmic
  Sobolev inequality with the Gaussian maximum of Shannon's entropy power.
  There is thus a complete parallel with the well-known link between
  logarithmic Sobolev inequalities and their information theoretic
  counterparts. We moreover provide an elementary proof of the reversed Gross
  inequality via a two-point inequality and the Central Limit Theorem.
\end{abstract}

\maketitle

\section{Shannon's entropy power and Gross's inequality}

In the sequel, we denote by $\ent{\mu}{f}$ the entropy of a
non-negative integrable function $f$ with respect to a positive measure $\mu$,
defined by
$$
\ent{\mu}{f} \defeq \int\!\!f\log f d\mu - \int\!\!f d\mu\,\log\int\!\!f d\mu.
$$
The Shannon
entropy \cite{shannon-1948} of an $n$-variate random vector $X$ with
probability density function (pdf) $f$ is given by
$$
 \sent{X} \defeq -\ent{\lan}{f}=-\int\!\!f\log f\,dx,
$$
where $dx$ denotes the $n$-dimensional Lebesgue measure on $\dR^n$. The
Shannon entropy power \cite{shannon-1948} of $X$ is then given by
$$
 \eent{X} \defeq \frac{1}{2\pi e}\exp\PAR{\frac{2}{n}\,\sent{X}}.
$$
It is well-known (cf. \cite{shannon-1948,cover-thomas-1991}) that Gaussians
saturates this entropy at fixed covariance. Namely, for any $n$-variate random
vector $X$ with covariance matrix $\K{X}$, one have
\begin{equation} \label{in:MNG}
 \eent{X} \leq \ABS{\K{X}}^{1/n},
\end{equation}
and $\ABS{\Kf}^{1/n}$ is the entropy power of the $n$-dimensional Gaussian
with covariance $\Kf$.

The logarithmic Sobolev inequality of Gross \cite{gross-1975} expresses that
for any non-negative smooth function $f : \dR^n \to \dR^+$
\begin{equation}
  \label{in:LSG}
  2\,\ent{\gn}{f} \leq \moy{\gn}{\frac{\ABS{\GR f}^2}{f}},
\end{equation}
where $\moyf{\gn}$ denotes the expectation with respect to $\gn$,
$\ABS{\cdot}$ the Euclidean norm and $\gn$ the $n$-dimensional standard
Gaussian given by
$$
 d\gn(x) \defeq (2\pi)^{-\frac{n}{2}}e^{-\frac{\ABS{x}^2}{2}}\,dx.
$$
Inequality~(\ref{in:LSG}) is sharp and the equality is achieved for $f$ of the
form $\exp(a \cdot)$.
 
By performing a change of function and an optimization, Beckner showed
\cite{beckner-1999} (see also \cite{carlen-1991}) that~(\ref{in:LSG}) is
equivalent to the following ``Euclidean'' logarithmic Sobolev inequality, for
any pdf $g$
\begin{equation} \label{in:SF}
  \ent{\lan}{g} 
   \leq 
  \frac{n}{2}\,\log\SBRA{\frac{1}{2\pi en}\int\!\!\frac{\ABS{\GR g}^2}{g}\,dx},
\end{equation}
where $\lan$ is the $n$-dimensional Lebesgue measure on $\dR^n$. Therefore,
for any $n$-variate random vector $X$ (with pdf $g$), we have
\begin{equation}\label{in:NJ}
 \eent{X}\fish{X} \geq n.
\end{equation}
This inequality can be obtained by many methods. The most classical ones are
via Shannon's entropy power inequality together with DeBruijn identity, or via
Stam's super-additivity of the Fisher information (cf.
\cite{stam-1959,carlen-1991,dembo-cover-thomas-1991,our-ls-2000}). Moreover,
Dembo showed in \cite{dembo-1990} that~(\ref{in:NJ}) is equivalent to
\begin{equation}\label{in:NJJ}
 \eent{X}\ABS{\fishm{X}}^{1/n} \geq 1, 
\end{equation}
where $\fishm{X}$ is the Fisher information matrix of $X$ defined by
$$
\fishm{X} \defeq \int\!\! \GR\log g \cdot \GR \log g^\top\,g\,dx,
$$
and we have $\fish{X}=\TR{\fishm{X}}$.  To deduce~(\ref{in:NJJ})
from~(\ref{in:NJ}), apply~(\ref{in:NJ}) to the random vector $X=\K{Y}^{-1/2}
Y$. Conversely, use the arithmetic-geometric means inequality
\begin{equation}\label{in:geomarith}
\PAR{a_1 \cdots a_n}^\frac{1}{n} \leq \frac{a+\cdots+a_n}{n}
\end{equation}
on the spectrum of the non-negative symmetric matrix $\fishm{X}$.

\section{Reversed Gross's logarithmic Sobolev inequality}

The Gross logarithmic Sobolev inequality~(\ref{in:LSG}) admits a reversed form
which states that for any positive smooth function $f : \dR^n \to \dR^+$
\begin{equation} \label{in:ILSG}
 \frac{\ABS{\moy{\gn}{\GR f}}^2}{\moy{\gn}{f}} \leq 2\,\ent{\gn}{f}.
\end{equation}
Here again, the $2$ constant is optimal and the equality is achieved for $f$
of the form $\exp(a \cdot)$.
Alike for~(\ref{in:LSG}), one can show by a change of function and an
optimization that the reverse form~(\ref{in:ILSG}) is equivalent to the
following inequality, for any pdf $g$
\begin{equation} \label{in:ISF}
 -\ent{\lan}{g} 
 \leq 
 \frac{n}{2} \log\SBRA{2\pi e \frac{\TR\K{g}}{n}},
\end{equation}
where $\K{g}$ is the covariance matrix of the pdf $g$.
Hence, we have for any $n$-variate random vector $X$ with pdf
\begin{equation}\label{in:NK}
 \eent{X} \leq \frac{\TR{\K{X}}}{n},
\end{equation}
where $\K{X}$ denotes the covariance matrix of $X$. This inequality is optimal
and is achieved by Gaussians $X$. Moreover, as we will show,
inequality~(\ref{in:NK}) is equivalent to~(\ref{in:MNG}).

Summarizing, we obtain the following statement
\begin{thm}
 The following assertions are true and equivalent
 \begin{enumerate}
 \item[(i)] For any smooth $f:\dR^n\to\dR^+$,
  $$
    \ABS{\moy{\gn}{\GR f}}^2 \leq 2\,\ent{\gn}{f}\moy{\gn}{f}.
  $$
  \item[(ii)] For any smooth $g:\dR^n\to\dR^+$, 
   $$
    -\ent{\lan}{g} \leq \frac{n}{2} \log\SBRA{\frac{2\pi e}{n}\TR\K{g}}.
   $$
 \item[(iii)] For any $n$-variate random vector $X$ with smooth pdf,
  $$
    n\eent{X} \leq \TR\K{X}.
  $$
 \item[(iv)] For any $n$-variate random vector $X$ with smooth pdf,
  $$
    \eent{X} \leq \ABS{\K{X}}^{1/n}.
  $$
  \end{enumerate}
\end{thm}
Therefore, there is a complete parallel between the equivalence
between~(\ref{in:LSG}), (\ref{in:SF}), (\ref{in:NJ}), (\ref{in:NJJ}) in one
hand and the equivalence between~(\ref{in:ILSG}), (\ref{in:ISF}),
(\ref{in:NK}), (\ref{in:MNG}) in the other hand.

\section{Sketches of proofs}

In this section, we present first two proofs of~(\ref{in:ILSG}), then we
explain how to deduce~(\ref{in:ISF}) from~(\ref{in:ILSG}) and (\ref{in:MNG})
from~(\ref{in:NK}) and vice versa.

The most natural way to establish~(\ref{in:ILSG}) is to start from a two-point.
inequality, just like Gross does for the logarithmic Sobolev
inequality~(\ref{in:LSG}) in \cite{gross-1975}.  Namely, if we denote by
$\beta$ the symmetric Bernoulli measure on $\BRA{-1,+1}$, one can show easily
that for any non-negative function $f:\BRA{-1,+1}\to {\dR^+}$,
\begin{equation}\label{eq:ils-bpq}
  (f(+1)-f(-1))^2 \leq \frac{1}{8}\; \ent{\beta}{f}\moy{\beta}{f}.
\end{equation}
This inequality is nothing else but the Csisz\'ar-Kullback inequality for
$\beta$ (see \cite{pinsker-1964}). Actually, the optimal
constant for the Bernoulli measure of parameter $p=1-q$ is $p^2q^2(\log q
-\log p)/(q-p)$, which resembles the optimal constant for the logarithmic
Sobolev inequality (see for example \cite{our-ls-2000}), but here again, only
the symmetric case gives the optimal constant of~(\ref{in:ILSG}).

The next step is to establish the following chain rule formula for $\entf{}$,
which generalizes the classical chain rule formula (cf.
\cite{shannon-1948,cover-thomas-1991}) for $\sentf$
\begin{prop}
  For any positive measures $\mu_i$ and their product $\mu$ on the product
  space, and for any bounded real valued measurable function $f$ on the
  product space, we have
\begin{equation}\label{in:entenso}
 \ent{\mu}{f} \geq \sum_{i=1}^n \ent{\mu_i}{\moy{\mu_{\bs i}}{f}},
\end{equation}
where $\mu_{\bs i}$ denotes the product of the measures $\mu_j$ with $j\neq i$.
\end{prop}

Finally, inequality~(\ref{in:ILSG}) can then be recovered by the use of the
Central Limit Theorem and integration by parts (in both discrete and Gaussian
forms). This concludes the first proof of~(\ref{in:ILSG}).

Actually, inequality~(\ref{in:ILSG}) can be recovered by a simple semi-group
argument, just like for the logarithmic Sobolev inequality (\ref{in:LSG}) (cf.
\cite{ledoux-zurich}). Namely, consider the heat semi-group
$\PAR{\SGf{t}}_{t\geq 0}$ on $\dR^n$, acting on a bounded continuous function
$f:\dR^n\to\dR$ as follows
$$
 \SG{t}{f}(x) 
  \defeq  \int_{\dR^n}\!\!f(x+\sqrt{t}\:y)\,d\gn(y).
$$
Notice that for any smooth function $f$, 
$\GR \SG{t}{f}=\SG{t}{\GR f}$ and
$$
 \partial_t\SG{t}{f}=\frac{1}{2}\:\LA \SG{t}{f}=\frac{1}{2}\:\SG{t}{\LA f}.
$$
Now, for any smooth positive bounded function $f:\dR^n\to \dR^+$, any $t\geq
0$ and any $x$, we can write, by performing an integration by parts and
omitting the $x$ variable
\begin{eqnarray*}
 \SG{t}{f\log f} - \SG{t}{f}\log\SG{t}{f} \
 &=&\ds{\int_0^t\!\!\partial_s\SBRA{\SG{s}{\SG{t-s}{f}\log\SG{t-s}{f}}}\,ds}\\
 &=&\ds{\frac{1}{2}\, 
        \int_0^t\!\!\SG{s}{\frac{\ABS{\GR\SG{t-s}{f}}^2}{\SG{t-s}{f}}}\,ds}.
\end{eqnarray*}
But by Cauchy-Schwarz inequality we get
$$
 \SG{s}{\frac{\ABS{\GR\SG{t-s}{f}}^2}{\SG{t-s}{f}}}
 =
 \SG{s}{\frac{\ABS{\SG{t-s}{\GR f}}^2}{\SG{t-s}{f}}}
 \geq
 \frac{\ABS{\SG{s}{\SG{t-s}{\GR f}}}^2}{\SG{s}{\SG{t-s}{f}}},
$$
which gives 
$$
 \SG{t}{f\log f} - \SG{t}{f}\log\SG{t}{f}
 \geq 
 \frac{t}{2} \frac{\ABS{\SG{t}{\GR f}}^2}{\SG{t}{f}}.
$$
Finally, inequality~(\ref{in:ILSG}) follows by taking $(t,x)=(1,0)$. Notice
that this method gives also the logarithmic Sobolev inequality (\ref{in:LSG}).
Namely, by Cauchy-Schwarz inequality
$$
 \ABS{\SG{t-s}{\GR f}}^2
 \leq
 \SG{t-s}{\ABS{\GR f}}^2
 \leq
 \SG{t-s}{f}\SG{t-s}{\frac{\ABS{\GR f}^2}{f}},
$$
therefore, we obtain
\begin{eqnarray*}
 \SG{t}{f\log f} - \SG{t}{f}\log\SG{t}{f} \
 &=&\ds{\frac{1}{2}\, 
    \int_0^t\!\!\SG{s}{\frac{\ABS{\SG{t-s}{\GR f}}^2}{\SG{t-s}{f}}}\,ds}\\ 
 &\leq&\ds{\frac{1}{2}\,
       \int_0^t\!\SG{s}{\SG{t-s}{\frac{\ABS{\GR f}^2}{f}}}\,ds}\\
 &=&\ds{\frac{t}{2}\,\SG{t}{\frac{\ABS{\GR f}^2}{f}}},
\end{eqnarray*}
which gives (\ref{in:LSG}) by taking here again $(t,x)=(1,0)$.

To deduce~(\ref{in:ISF}) from~(\ref{in:ILSG}), just apply~(\ref{in:ILSG}) to
$$
 f(x)=h(x)\,\PAR{2\pi}^\frac{n}{2}\,e^{\frac{\ABS{x}^2}{2}} 
$$
where $h$ is a compactly supported smooth pdf. One then gets
$$
 \ABS{\int\!\!x\,h\,dx+\int\!\!\GR h\,dx}^2
 \leq
 2\!\int\!\!h\log h\,dx + \int\!\!\ABS{x}^2\,h\,dx 
 +n\log(2\pi).
$$
But we have $\int\!\!\GR h\,dx=0$. Therefore, by denoting $\K{h}$ the
covariance matrix associated with the pdf $h$, one gets
$$
 -\ent{\lan}{h} 
 \leq 
 \frac{1}{2} \TR\K{h} + \frac{n}{2}\log(2\pi),
$$
which remains true for any smooth pdf $h$. Finally, by performing the change
of function $h=\alpha g(\alpha \cdot)$ and optimizing in $\alpha$, one obtains

$$
 -\ent{\lan}{h} \leq \frac{n}{2}\log\SBRA{2\pi e\frac{\TR\K{h}}{n}}
$$
which is nothing else than~(\ref{in:NK}). Conversely, it is easy to see that
we can recover~(\ref{in:ILSG}) for any pdf $f$ by approximating $f$ by
compactly supported probability density functions.

The equivalence between~(\ref{in:NK}) and~(\ref{in:MNG}) is obtained as for
the equivalence between~(\ref{in:NJ}) and~(\ref{in:NJJ}). Namely, to
deduce~(\ref{in:MNG}) from~(\ref{in:NK}), apply~(\ref{in:NK}) to the random
vector $X=\K{Y}^{-1/2} Y$. Conversely, use the arithmetic-geometric means
inequality (\ref{in:geomarith}) on the spectrum of the non-negative symmetric
matrix $\K{X}$.

\section{Remarks}

It is well-know that the logarithmic Sobolev inequality~(\ref{in:LSG}) is a
consequence of the Gaussian isoperimetric inequality \cite{ledoux-zurich}. In
contrast, it is shown in \cite{barthe-cordero-fradelizi} that the reversed
form~(\ref{in:ILSG}) is equivalent to a translation property
\cite{bobkov-transl}. Namely, for any smooth function $f : \dR^n \to [0,1]$
\begin{equation}\label{in:bobkov} 
 \ABS{\moy{\gn}{\GR f}} \leq \isop{\moy{\gn}{f}},
\end{equation}
where $\isopf$ is the Gaussian isoperimetric function given by $\isopf\defeq
\Phi' \circ \Phi^{-1}$, where $\Phi$ is the Gaussian distribution function
given by $\Phi(\cdot)\defeq \g((-\infty,\cdot])$. Bobkov's
inequality~(\ref{in:bobkov}) expresses that among all measurable sets with
fixed Gaussian measure, half spaces have minimum barycenter
\cite{barthe-cordero-fradelizi}.

\subsubsection*{Acknowledgments}

The author would like to thanks Michel Ledoux and C\'ecile An\'e for their
encouragements and helpful comments.


\begin{thebibliography}{10}

\bibitem{our-ls-2000}
C.~An\'e, S.~Blach\`ere, D.~Chafa{\"{\i}}, P.~Foug\`eres, I.~Gentil,
  F.~Malrieu, C.~Roberto, and G.~Scheffer.
\newblock Sur les in\'egalit\'es de {S}obolev logarithmiques.
\newblock to appear in ``Panoramas et Synth\`eses, Soci\'et\'e Math\'ematique
  de France, 2001.

\bibitem{bakry-concordet-ledoux-1997}
D.~Bakry, D.~Concordet, and M.~Ledoux.
\newblock Optimal heat kernel bounds under logarithmic {S}obolev inequalities.
\newblock {\em ESAIM Probab. Statist.}, 1:391--407 (electronic), 1997.

\bibitem{barthe-cordero-fradelizi}
F.~Barthe, D.~Cordero-Erausquin, and M.~Fradelizi.
\newblock Shift inequalities of gaussian type and norms of barycenters.
\newblock preprint, september 1999.

\bibitem{beckner-1999}
W.~Beckner.
\newblock Geometric asymptotics and the logarithmic {S}obolev inequality.
\newblock {\em Forum Math.}, 11(1):105--137, 1999.

\bibitem{blachman-1965}
N.~Blachman.
\newblock The convolution inequality for entropy powers.
\newblock {\em IEEE Trans. Information Theory}, IT-11:267--271, 1965.

\bibitem{bobkov-transl}
S.~Bobkov.
\newblock The size of singular component and shift inequalities.
\newblock {\em Ann. Probab.}, 27:416--431, 1999.

\bibitem{carlen-1991}
E.~Carlen.
\newblock Super-additivity of {F}isher's information and logarithmic {S}obolev
  inequalities.
\newblock {\em J. Funct. Anal.}, 101(1):194--211, 1991.

\bibitem{cover-thomas-1991}
T.~Cover and J.~Thomas.
\newblock {\em Elements of information theory}.
\newblock John Wiley \& Sons Inc., New York, 1991.
\newblock A Wiley-Interscience Publication.

\bibitem{dembo-1990}
A.~Dembo.
\newblock Information inequalities and uncertainty principles.
\newblock In {\em Tech. Rep., Dept. of Statist.} Stanford Univ., 1990.

\bibitem{dembo-cover-thomas-1991}
A.~Dembo, T.~Cover, and J.~Thomas.
\newblock Information-theoretic inequalities.
\newblock {\em IEEE Trans. Inform. Theory}, 37(6):1501--1518, 1991.

\bibitem{gross-1975}
L.~Gross.
\newblock Logarithmic {S}obolev inequalities.
\newblock {\em Amer. J. Math.}, 97(4):1061--1083, 1975.

\bibitem{ledoux-berlin}
M.~Ledoux.
\newblock Concentration of measure and logarithmic {S}obolev inequalities.
\newblock In {\em S\'eminaire de Probabilit\'es, XXXIII}, Lecture Notes in
  Math., pages 120--216. Springer, Berlin, 1999.

\bibitem{ledoux-zurich}
M.~Ledoux.
\newblock The geometry of {M}arkov {D}iffusion {G}enerators.
\newblock {\em Ann. Fac. Sci. Toulouse Math.}, IX(2):305--366, 2000.

\bibitem{pinsker-1964}
M.~S. Pinsker.
\newblock {\em Information and information stability of random variables and
  processes}.
\newblock Holden-Day Inc., San Francisco, Calif., 1964.
\newblock Translated by Amiel Feinstein.

\bibitem{shannon-1948}
C.~Shannon.
\newblock A mathematical theory of communication.
\newblock {\em Bell System Tech. J.}, 27:379--423, 623--656, 1948.

\bibitem{stam-1959}
A.~Stam.
\newblock Some inequalities satisfied by the quantities of information of
  {F}isher and {S}hannon.
\newblock {\em Information and Control}, 2:101--112, 1959.

\end{thebibliography}
\end{document}